\documentclass[12pt,leqno]{article}
\usepackage[francais]{babel}
\usepackage{amsmath,amsfonts,amssymb,amsthm}
\begin{document}

\newcommand\1{\mbox{1\hspace{-.3em}I}}
\newcommand\Pp{\mathbf{P}}
\newcommand\Z{\mathbb{Z}}
\newcommand\A{\mathbb{A}}
\newcommand\F{\mathbb{F}}
\newcommand\C{\mathbb{C}}
\newcommand\N{\mathbb{N}}
\newcommand\R{\mathbb{R}}
\newcommand\B{\mathbb{B}}
\newcommand\G{\mathbb{G}}
\newcommand\Gg{\textbf{G}}
\newcommand\Q{\mathbb{Q}}
\newcommand\Ql{\mathbb{Q}_\ell}
\newcommand\oql{\overline{\mathbb{Q}}_\ell}
\newcommand\Lc{\Lambda^\bullet}
\newcommand\T{\mathbb{T}}
\newcommand\Hh{\mathfrak{H}}
\newcommand\dd{\mathfrak{d}}
\newcommand\Ss{\mathcal{S}}
\newcommand\bb{\mathcal{B}}
\newcommand\Nn{\mathcal{N}}
\newcommand\Dd{\mathcal{D}}
\newcommand\Ee{\mathcal{E}}
\newcommand\Ff{\mathcal{F}}
\newcommand\Ll{\mathcal{L}}
\newcommand\Oo{\mathcal{O}}
\newcommand\alg{\mathcal{A}lg}
\newcommand\Bb{\mathcal{B}}
\newcommand\Uu{\mathcal{U}}
\newcommand\Qq{\mathcal{Q}}
\newcommand\ind{\mathrm{ind}}
\newcommand\Frob{\mathrm{Frob}}
\newcommand\rg{\rightarrow}
\newcommand\lgr{\longrightarrow}
\newcommand\Gal{\mathrm{Gal}}
\newcommand\End{\mathrm{End}}
\newcommand\Res{\mathrm{Res}}
\newcommand\Ind{\mathrm{Ind}}
\def\adots{\mathinner{\mkern1mu\raise1pt\vbox{\kern7pt\hbox{.}}
\mkern2mu\raise4pt\hbox{.}
\mkern2mu\raise7pt\hbox{.}\mkern1mu}}

\title{Jacques Peyri\`ere et les produits de Riesz}
\author{Jean--Pierre Kahane}
\date{}
\maketitle

The present article,\og Jacques Peyri\`ere et les produits de Riesz\fg, was supposed to be a chapter of a book devoted to Jacques Peyri\`ere, to be published by Birkha\"user. I was asked to sign a CAP form
(contributor's autorization to publish) and was willing to give such an autorization, but the form meant that I was giving up all rights on my work, including translation. I have this paper written in French on purpose and I don't want to see it translated in English without my permission. I asked for a change, but no change was possible, and all contributors had to agree with the following statement :
\og The contributor hereby irrevocably grants and assigns to Springer any and all rights which the contributor may have, or may at any time be found to have, in and to the Article, including but not limited to the sole right to print, publish and sell the Article throughout the world, and including all rights in and to all revisions or versions in all languages and media throughout the world.\fg\ Since I refused to sign the CAP as it was proposed, the publisher decided, not to publish the paper.
I am sorry, of course, not to participate directly to the homage paid to Jacques Peyri\`ere. As a personal homage I should add a few words.
The copyright of scientific papers should belong to the authors or to scientific bodies. When a publisher is willing to publish a paper of ours, it is natural that he asks us for permission, and natural that we grant him this permission, with thanks if the publication is well done. But when the publisher insists to have all rights on the work, he converts the copyright into a commercial good ; the author can discover that a copyright granted to a publisher at a given time belongs now to another publisher, and that it was sold indeed by the first publisher as part of his \og fonds de commerce.\fg\ 
The publishers of scientific books or journals played an important and positive role in the past, and may be now as useful than before, in a different context. But being too eager to grasp everything from the authors may be wrong at a long run not only for the authors but for themselves. The context of open archives may create a new need for regular publications, if good relations are created between authors and publishers. If the most powerful publishers do not fill this condition, they may have difficulties in the future with the scientific communities. 

\vskip6mm
L'\oe uvre de Jacques Peyri\`ere s'\'etend sur des domaines vari\'es de l'analyse harmonique, de la th\'eorie g\'eom\'etrique de la mesure, des probabilit\'es, de la combinatoire et de l'alg\`ebre. Il y a des traits communs \`a tous ces travaux : l'audace de s'attaquer \`a des probl\`emes difficiles, l'\'el\'egance de la pr\'esentation et du style, et l'influence  \`a long terme des m\'ethodes qu'il a introduites. Tous ces \'el\'ements sont pr\'esents d\`es le premier grand article de Jacques Peyri\`ere, sur les produits de Riesz \cite{pey1}. On pourrait \'ecrire un livre sur tout ce qui en a d\'ecoul\'e. Je me bornerai \`a en d\'egager quelques traits saillants et quelques prolongements.

Les produits de Riesz que consid\`ere Peyri\`ere sont de la forme
\begin{equation}
\prod_{j=0}^\infty (1+\textrm{Re}(a_j\ e^{i\lambda_j t}))
\end{equation}
o\`u $t$ est la variable r\'eelle, des $\lambda_j$ des entiers $\ge 1$ tels que
\begin{equation}
\frac{\lambda_{j+1}}{\lambda_j} \ge 3\qquad (j=0,1,2,\ldots)
\end{equation}
et les $a_j$ des nombres complexes de modules $\le 1$
\begin{equation}
|a_j| \le 1 \qquad (j=0,1,2,\ldots)\,.
\end{equation}
La condition (3) signifie que tous les facteurs sont $\ge 0$, et la condition (2) entra\^{i}ne que dans le d\'eveloppement  du produit (1) en somme il n'y a pas d'interf\'erence entre les fr\'equences. Ainsi les produits partiels de~(1)
\begin{equation}
P_{a,N}(t) = \prod_0^N (1+\textrm{Re} (a_j e^{i\lambda_j t}))
\end{equation}
s'\'ecrivent, en posant $a_j=r_j e^{i\theta_j}$,
\begin{equation}
\sum_{(\varepsilon_i)}\big(\frac{r_0}{2}\big)^{|\varepsilon_0|}  \big(\frac{r_1}{2}\big)^{|\varepsilon_1|} \cdots \big(\frac{r_N}{2}\big)^{|\varepsilon_N|}\ e^{i(\varepsilon_0\theta_0+\varepsilon_1\theta_1+\cdots+\varepsilon_N\theta_N)} \ e^{i(\varepsilon_0\lambda_0+\varepsilon_1\lambda_1+\cdots+\varepsilon_N\lambda_N)t}
\end{equation}
o\`u les $\varepsilon_j$ prennent les valeurs $-1,\ 0,\ 1$, et les $\varepsilon_0\lambda_0+\varepsilon_1\lambda_1+\cdots+\varepsilon_N\lambda_N$ sont tous distincts. (5) est donc la somme partielle d'ordre $\lambda_0+\lambda_1+\cdots+\lambda_N$ de la s\'erie
\begin{equation}
\sum_{(\varepsilon_j)\in\{-1,0,1\}^\N,\sum|\varepsilon_j|<\infty} \prod\big(\frac{r_j}{2}\big)^{|\varepsilon_j|} e^{i\sum\varepsilon_j\theta_j} e^{i (\sum \varepsilon_j\lambda_j)t}
\end{equation}
Comme ces sommes partielles sont positives et de valeur moyenne 1, (6) est la s\'erie de Fourier d'une mesure qu'on d\'esigne par $\mu_a$, et on identifie le produit de Riesz \`a $\mu_a$.

On consid\`ere les $\lambda_j$ comme fix\'es, par exemple $\lambda_j=4^j$ (c'\'etait le choix initial de Fr\'ed\'eric Riesz) et on \'etudie les propri\'et\'es de $\mu_a$ en fonction de $a=(a_j)$ sous la condition (3).

D'abord, en appliquant le crit\`ere de Wiener \`a la s\'erie (6), on voit que $\mu_a$ est une mesure diffuse. Voici un \'echantillon des r\'esultats de Peyri\`ere :

\vskip2mm

A) (th\'eor\`eme 1.2) Supposons $|a_j|\le 1$, $|b_j|\le 1$ $(j=0,1,2,\ldots)$ et
\begin{equation}
\sum_0^\infty |a_j-b_j|^2 =\infty\,.
\end{equation}
Alors les deux mesures $\mu_a$ et $\mu_b$ sont \'etrang\`eres : $\mu_a \bot \mu_b$.

\vskip2mm

B) (cas particulier du th\'eor\`eme 1.3). Supposons $\lambda_j=4^j$, $|a_j|\le 1$, $|b_j\le 1$ $(j=0,1,2,\ldots)$, et une condition plus forte que
\begin{equation}
\sum_0^\infty |a_j -b_j|^2 <\infty\,.
\end{equation}
\`a savoir
$$
\sum_0^\infty   |a_j-b_j | ^2/ (1-|a_j|)<\infty\,.
\leqno(8')
$$
Alors $\mu_a$ et $\mu_b$ sont mutuellement absolument continues : $\mu_a\sim \mu_b$.

\vskip2mm

C) (cas particulier du th\'eor\`eme 2.8). Supposons encore $\lambda_j=4^j$, et $|a_j|\le 1$ $(j=0,1,2,\ldots)$. Alors les dimensions de Hausdorff des bor\'eliens $E$ que charge $\mu_a$ $(\mu_a(E)>0)$ v\'erifient
\begin{eqnarray}
1 &-\limsup\limits_{n\rg \infty}\frac{1}{\log \lambda_n} \int \log P_{a,n}(t) d\mu_a(t) \le \dim E\\
\nonumber&\le 1-\liminf\limits_{n\rg \infty} \frac{1}{\log \lambda_n} \int \log P_{a,n}(t) d\mu_a(t)
\end{eqnarray}

D) La d\'efinition de $\mu_a$ s'\'etend au cas $\lambda_j=2^j$, sous la condition additionnelle
\begin{equation}
\sup |a_j| <1\,,
\end{equation}
et l'\'enonc\'e A est valable sous la condition additionnelle $\sup|a_j|<1$,\break $\sup|b_j|<1$.

\vskip2mm

Preuve de A et commentaires.

Voici la preuve de A par Peyri\`ere. On v\'erifie que les fonctions
\begin{equation}
e^{i\lambda_j t}- \frac{1}{2}\  \overline{a}_j \qquad (j=0,1,2,\ldots)
\end{equation}
forment un syst\`eme orthonormal dans $L^2(\mu_a)$, et que les normes sont comprises entre deux nombres strictement positifs. On choisit une suite $\{c_j\}\in \ell^2$ telle que
\begin{eqnarray}
\left\{\begin{array}{l}
c_j(\overline{a}_j-\overline{b}_j) \ge 0   \\
\noalign{\vskip2mm}
\sum c_j(\overline{a}_j-\overline{b}_j) = \infty  \\
\end{array}
\right.
\end{eqnarray}
ce qui est possible gr\^ace \`a (7), et on consid\`ere les deux s\'eries
\begin{equation}
\sum_0^° c_j(e^{i\lambda_j t} - \frac{1}{2}\ \overline{a}_j)\,,
\end{equation}
\begin{equation}
\sum_0^° c_j(e^{i\lambda_j t} - \frac{1}{2}\ \overline{b}_j)\,.
\end{equation}
La premi\`ere converge dans $L^2(\mu_a)$ et la seconde dans $L^2(\mu_b)$. Il existe donc une suite d'entiers $N$ tendant vers l'infini telle que les sommes partielles d'ordre $N$ de (13) convergent $\mu_a-p\cdot p\cdot$ et les sommes partielles d'ordre $N$ de (14) convergent $\mu_b -p\cdot p\cdot$. Or, d'apr\`es (12), elles ne peuvent pas converger en un m\^eme point $t$. Donc $\mu_a$ et $\mu_b$ sont \'etrang\`eres.

Presque en m\^eme temps que Peyri\`ere et ind\'ependamment, G. Brown et W. Moran avaient obtenu l'\'enonc\'e A \cite{brmo}. L'\'el\'egance de la preuve de Peyri\`ere est frappante.

On peut adapter cette preuve \`a des situations un peu diff\'erentes, o\`u apparaissent des produits de Riesz g\'en\'eralis\'es. Par exemple, c'est un bon outil pour comprendre la construction faite par Katznelson d'une s\'erie trigonom\'etrique dont toutes les sommes partielles sont positives (ce qui entra\^{i}ne, par un th\'eor\`eme de Helson, que les coefficients tendent vers 0) et qui n'est pas une s\'erie de Fourier--Lebesgue \cite{kat}. Je l'ai utilis\'e explicitement pour la construction d'une s\'erie de Walsh ayant la m\^eme propri\'et\'e~\cite{kah}.

La preuve de Peyri\`ere repose sur le fait que la s\'erie (13) converge dans $L^2(\mu_a)$ lorsque $\{c_j\}\in \ell^2$ (on peut d'ailleurs la pr\'esenter en disant que, sous l'hypoth\`ese $\{c_j\}\in \ell^2$, (13) converge dans $L^2(\mu_a)$ et (14) dans $L^2(\mu_b)$, donc (13) et (14) dans $L^2(\mu_a\wedge \mu_b)$, donc si $\mu_a\wedge \mu_b\neq0$, $\sum c_j(\overline{a}_j - \overline{b}_j)$ est une s\'erie convergente d\`es que $\{c_j\}\in \ell^2$, ce qui entra\^{i}ne $\{a_b-b_j\}\in \ell^2$). On peut se demander si la s\'erie (13) converge $\mu_a$--presque partout. Il en est bien ainsi, comme l'ont montr\'e ind\'ependamment J.~Peyri\`ere \cite{pey2} et A.H.~Fan \cite{fan1}. La preuve de Fan utilise une randomisation du produit de Riesz, \`a savoir une mesure al\'eatoire
\begin{equation}
\mu_{a,\omega} \sim \prod_{j=0}^\infty (1+ \textrm{Re}(a_j\ e^{i\omega_j}\ e^{i\lambda_j t}))
\end{equation}
o\`u les $\omega_j$ sont des variables al\'eatoires ind\'ependantes distribu\'ees sur $[0,2\pi]$ suivant la mesure de Lebesgue normalis\'ee.

\vskip2mm

Commentaires sur B.

La conclusion $\mu_a\sim \mu_b$ est tr\`es facile \`a obtenir \`a partir de (8), pour toute suite $(\lambda_j)$ v\'erifiant (2), si l'on fait l'hypoth\`ese (10). Mais, sans cette hupoth\`ese, les conditions (2), (3) et (8) ne garantissent pas $y_a\sim y_b$. Le th\'eor\`eme~1.3 de Peyri\`ere donne d'autres conditions suffisantes. F.~Parreau a montr\'e en 1989~que
\begin{equation}
|a_j| = |b_j|\qquad (j=0,1,2,\ldots)
\end{equation}
joint \`a (2), (3) et (8) entra\^{i}ne $\mu_a\sim \mu_b$ \cite{par}. 

Voici un th\'eor\`eme de Kilmer et Saeki de 1988, relatif aux produits de Riesz al\'eatoires (15). Si l'on munit le disque $\{z=\textrm{re}^{i\theta},\ 0\le r\le 1\}$ de la m\'etrique
\begin{equation}
ds^2 =d\theta^2 +(1-r)^{-1/2} dr^2\,,
\end{equation}
la condition
\begin{equation}
\sum_0^\infty d^2 (a_j,b_j)<\infty
\end{equation}
(plus forte que (8)), jointe \`a (2) et (3), entra\^{i}ne que $\mu_{a,\omega}\sim \mu_{b,\omega}$ presque s\^urement En corollaire, si on fait au lieu de (2), l'hypoth\`ese plus forte 
\begin{equation}
\sum_0^\infty \Big(\frac{\lambda_{j+1}}{\lambda_j}\Big)^2 < \infty \,,
\end{equation}
(18) entra\^{i}ne $\mu_a\sim \mu_b$
\cite{kisa}.

Fan a pos\'e la question si
\begin{equation}
\mu_{a,\omega} \sim \mu_{b,\omega}\  p\cdot s\cdot \Longleftrightarrow \mu_a \sim \mu_b \ \cite{fan2}\,.
\end{equation}
Si la r\'eponse \'etait positive, (18), qui contient la condition de Parreau (16), serait la meilleure condition pour entra\^{i}ner $\mu_a\sim \mu_b$ sous les conditions (2)~et~(3).

\vskip2mm

Commentaire sur C.

Peyri\`ere a \'et\'e un virtuose dans la d\'etermination des dimensions de Hansdorff des mesures positives. Une mesure $\mu$ \'etant donn\'ee, il est naturel d'introduire une dimension inf\'erieure et une dimension sup\'erieure, qui encadrent les dimensions des bor\'eliens que charge la mesure. Lorsqu'elles sont \'egales, on parle de dimension de la mesure.

Par exemple, C montre que, dans l'hypoth\`ese $\lambda_j=4^j$, si l'on suppose l'existence de la limite
\begin{equation}
L=\lim_{n\rg \infty} \frac{1}{\log \lambda_n} \int \ \log \ P_{a,n}(t)\ d\mu_a(t)\,,
\end{equation}
$\mu_a$ a pour dimension $1-L$.

Le th\'eor\`eme 2.8 de Peyri\`ere d\'epasse ce cas particulier : au lieu de $\lambda_j=4^j$, il suffit de supposer que les quotients $\frac{\lambda_{j+1}}{\lambda_j}$ sont entiers, $\ge 3$, et born\'es.

La m\'ethode de Peyri\`ere repose sur la th\'eorie ergodique d\'evelopp\'ee par Billingsley. C'est aussi de cette mani\`ere qu'il a calcul\'e la dimension des mesures al\'eatoires obtenues \`a partir des cascades multiplicatives de Beno\^{i}t Mandelbrot \cite{kape}. On doit s'attendre, et c'est le cas, \`a des r\'esultats plus simples \`a exprimer pour des mesures al\'eatoires que pour des mesures d\'eterministes.

Les membres extr\^emes de la formule (9) donnent des estimations de $\dim_*\mu_a$ et $\dim^*\mu_a$, les dimensions inf\'erieure et sup\'erieure de $\mu_a$, et le th\'eo\-r\`e\-me~2.8 de Peyri\`ere est le premier exemple de ce type d'estim\'e. L'analyse dimensionnelle des mesures (non n\'ecessairement donn\'ees par un produit de Riesz) consiste \`a mettre en \'evidence, pour chaque $\alpha$ compris entre $\dim_*\mu$ et $\dim^*\mu$, un bor\'elien $B_\alpha$ de dimension $\alpha$ tel que $\mu(B_\alpha)$ soit le plus grand possible, de fa\c con que $\alpha<\alpha'$ entra\^{i}ne $B_\alpha \subset B_{\alpha'}$. Le chapitre 5 de l'expos\'e de Barral, Fan et Peyri\`ere \`a para\^{i}tre \`a Panoramas et Synth\`eses \cite{bafape} fait le point sur l'analyse dimensionnelle. Un bon outil est l'exposant de H\"older local
\begin{equation}
\alpha(t) = \liminf_{h\rg0}\ \frac{\log \mu([t-h,t+h])}{\log h}
\end{equation}
que Fan appelle la dimension locale inf\'erieure de $\mu$ au point $t$. Les formules de~Fan
\begin{eqnarray}
\dim_*\mu &= \sup \{\alpha\ge 0 : \alpha(t) \ge \alpha\ \ \mu-pp\}\\
\dim^*\mu &= \inf \{\alpha\ge 0 : \alpha(t) \le \alpha\ \ \mu-pp\}
\end{eqnarray}
\'etablissent la relation du local au global.

L'analyse multifractale, \`a laquelle Jacques Peyri\`ere a beaucoup contribu\'e, repose sur la consid\'eration de la fonction $\alpha(t)$ (souvent restreinte dans la litt\'erature aux $t$ tels que dans (22) la limite inf\'erieure co\"{i}ncide avec la limite sup\'erieure). Elle consiste \`a \'etudier la dimension $D(\alpha)$ des ensembles
\begin{equation}
E_\alpha = \{t : \alpha(t) = \alpha\}\,.
\end{equation}
\cite{bafape} contient une \'etude approfondie de l'analyse multifractale, mais dit peu de choses sur l'analyse multifractale des produits de Riesz, qui semble un sujet difficile.

Le d\'eveloppement du produit (1) en s\'erie trigonom\'etrique (6) permet des approches, assez grossi\`eres, que je vais indiquer maintenant.

\vskip2mm

Premi\`ere approche. L'examen de la s\'erie (6) donne tout de suite $\widehat{\mu}_a \widehat{\mu}_b=\widehat{\mu}_{\frac{ab}{2}}$, c'est--\`a--dire $\mu_a \ast \mu_b =\mu_{\frac{ab}{2}}$. Si $\mu_a$ charge un bor\'elien $A$ et $\mu_b$ un bor\'elien $B$, $\mu_{\frac{ab}{2}}$ charge donc la somme $A+B$. Comme $\dim(A+B) \le \dim A +\dim B$, on a
\begin{equation}
\dim_* \mu_{\frac{ab}{2}} \le \dim_* \mu_a+ \dim_*\mu_b\,.
\end{equation}
De la m\^eme fa\c con
\begin{equation}
\dim^* \mu_{\frac{ab}{2}} \le \dim^* \mu_a +\dim^* \mu_b\,.
\end{equation}
Voici une application imm\'ediate. Supposons $\lambda_j=4^j$ $(j=0,1,2,\ldots)$. Alors le spectre de $\mu_a$ est contenu dans la r\'eunion des intervalles ouverts $]\frac{2}{3} A^n,\ \frac{4}{3} 4^n[$. L'image $\mu_b$ de $\mu_a$ par la transformation $t\lgr 2t$ a son spectre dans la r\'eunion des intervalles $]\frac{1}{3} 4^n,\ \frac{2}{3} 4^n[$. Comme ces spectres sont disjoints, $\mu_a\ast \mu_b$ est la mesure de Lebesgue, dont la dimension est 1, et comme les dimensions de $\mu_a$ et de $\mu_b$ sont les m\^emes, on~a
\begin{equation}
\dim_* \mu_a \ge \frac{1}{2}\,.
\end{equation}
Cela vaut quels que soient les $a_s$.

\vskip2mm

Seconde approche. L'\'energie d'une mesure $\mu$ d\'efinie sur le cercle par rapport au noyau de potentiel $|tg \frac{t-s}{2}|^{-\alpha}$ ($\alpha$--\'energie)~est
\begin{equation*}
\Ee_\alpha(\mu) = \int\!\int\frac{d\mu(t)d\mu(s)}{|tg\frac{t-s}{2}|^\alpha}\quad (0<\alpha<1)\,.
\end{equation*}
Ainsi $\mu$ est d'$\alpha$--\'energie finie si et seulement si
\begin{equation}
\sum_{n\in \Z} |\widehat{\mu}(n)|^2 |n|^{\alpha-1} < \infty\,.
\end{equation}
Sous l'hypoth\`ese
\begin{equation}
\sup \frac{\lambda_{j+1}}{\lambda_j} < \infty\ ,
\end{equation}
$\mu_a$ est donc d'$\alpha$--\'energie finie si et seulement si
\begin{equation*}
\sum_n \Big(\prod_{0}^{n+1} (1+|a_j|^2) - \prod_0^n (1+|a_j|^2))\ \lambda_n^{\alpha-1} < \infty\,,
\end{equation*}
soit
\begin{equation}
\sum_n \lambda_n^{\alpha-1} \prod_0^n (1+|a_j|^2) < \infty\,.
\end{equation}
Pour fixer les id\'ees, prenons $\lambda_j=4^j$ et $a_j=a>0$. Alors la condition pour que $\mu_a$ soit d'$\alpha$--\'energie finie est
\begin{equation}
4^{\alpha-1}(1+a^2) < 1\,.
\end{equation}
On sait qu'un bor\'elien est d'$\alpha$--capacit\'e positive si et seulement s'il porte une mesure d'$\alpha$--\'energie finie. Si $\mu_a$ charge un bor\'elien, ce bor\'elien est donc d'$\alpha$--capacit\'e positive lorsque $\alpha$ v\'erifie (32), et cela entra\^{i}ne que sa dimension est $\ge 1 -\frac{\log (1+a^2)}{\log 4}$. Donc
\begin{equation}
\dim_* \mu_a \ge 1 - \frac{\log(1+a^2)}{\log 4}\,,
\end{equation}
ce qui am\'eliore (28) mais n'atteint pas la valeur donn\'ee par (21).

\vskip2mm

Troisi\`eme approche. On peut chercher \`a lire sur la s\'erie (6) des indications sur la fonction $\alpha(t)$ d\'efinie en (22). Pour fixer les id\'ees prenons $\lambda_j=4^j$, et posons $\nu_n=\lambda_n -\lambda_{n-1}-\lambda_{n-2}-\cdots -\lambda_0$. Ainsi l'intervalle $]\nu_n,2\nu_n]$ est disjoint du spectre de $\mu_a$. Soit $V_n$ le noyau de de la Vall\'ee Poussin d'ordre $2\nu_n$ : ses coefficients valent 1 sur $[-\nu_n,\nu_n]$, et 0 hors de $[-2\nu_n,2\nu_n]$, et c'est une combinaison lin\'eaire \`a coefficients $-1$ et 2 des noyaux de Fej\'er d'ordres $\nu_n$ et $2\nu_n$. On a 
\begin{eqnarray}
P_{n,a}(t)   &=& \int_{-\pi}^\pi V_n(s)d\mu_a(t-s)\\
\noalign{\vskip2mm}
\nonumber &=& \int_0^\pi V_n' (s) \mu_a([t-s,t+s]) ds\\
\noalign{\vskip2mm}
\nonumber &\le& 10 \int_0^\pi \big(\frac{\nu_n}{s} \wedge \frac{1}{s^2}\Big) \mu_a([t-s,t+s]) ds\,,
\end{eqnarray}
donc on a l'implication
\begin{equation}
\forall s\ \mu_a([t-s,t+s]) \le C\ s^{\beta(t)} \Longrightarrow \forall n \ P_{a,n}(t) \le C' \lambda_n^{1-\beta(t)}\,,
\end{equation}
$\frac{C'}{C}$ ne d\'ependant que de la suite $(\lambda_j)$ ; au lieu de $\lambda_j=4^j$ on peut prendre n'importe quelle suite $(\lambda_j)$ v\'erifiant
\begin{equation}
\inf \frac{\lambda_{j+1}}{\lambda_j} > 3\,,\ \ \sup  \frac{\lambda_{j+1}}{\lambda_j}  < \infty\,.
\end{equation}

Dans l'autre sens les calculs sont moins simples et les r\'esultats moins bons. Partons des hypoth\`eses (2) et (3), et d\'efinissons
\begin{equation}
L_\nu(t)=\frac{1}{\nu} + 2 \sum_{n=1}^\infty \frac{\cos nt}{\nu+n}\,.
\end{equation}
C'est une fonction $\ge 0$, de valeur moyenne $\frac{1}{\nu}$. Posons ici
\begin{equation}
\nu_j=\lambda_{j+1} - \lambda_j - \lambda_{j-1} - \cdots - \lambda_0\,;
\end{equation}
ainsi les fr\'equences de $P_{a,j+1} - P_{a,j}$ sont $\ge \nu_j$. D'apr\`es (6),
\begin{eqnarray*}
 {\mu([t-s,t+s]) = \int_{t-s}^{t+s} P_{a,N}(u)du}\hskip10mm \\
+ \sum_{j\ge N} \sum_n \widehat{P_{a,j+1}-P_{a,j}}(n) \frac{e^{in(t+s)} -e^{in(t-s)}}{in}\,.
\end{eqnarray*}
Posons
\begin{equation}
\Qq_{a,j} = (P_{a,j+1}+P_{a,j}) \ast L_{\nu_j}\,.
\end{equation}
Alors
\begin{equation}
\mu([t-s,t+s]) \le \int_{t-s}^{t+s} P_{a,N}(u)du+\sum_{j\ge N}(\Qq_{a,j}(t+s) +\Qq_{a,j}(t-s))
\end{equation}

Voici une application simple de la formule (40). Supposons (30) en plus de (2) et (3), et
\begin{equation}
\forall t \forall j P_{a,j}(t) \le C\ \lambda_j^{1-\beta}\, ;
\end{equation}
alors
\begin{equation}
\forall t \forall s \ \mu([t-s,t+s]) \le C' s^{-\beta}
\end{equation}
o\`u $\frac{C'}{C}$ ne d\'epend que de la suite $(\lambda_j)$.

\vskip 2mm

Commentaire sur D, et digression

La suite $2^j$ a la propri\'et\'e qu'une somme finie $\sum \varepsilon_j 2^j$ $(\varepsilon_j \in \{-1,0,1\})$ est nulle si et seulement si tous les $\varepsilon_j$ sont nuls. A partir de cette observation Peyri\`ere montre qu'en prenant $\lambda_j =2^j$ les produits partiels (4) sont born\'es dans $L^1$ et que, sous la condition (10), ils convergent vaguement vers une mesure $\mu_a$, qu'on peut encore identifier au produit de Riesz~(1).

Suivant cette id\'ee, un cadre naturel pour l'\'etude des produits de Riesz est, au lieu de la condition (3), l'hypoth\`ese de quasi--ind\'ependance des $\lambda_j$, \`a savoir que
\begin{equation}
\sum_{\textrm{finie}} \varepsilon_j\lambda_j =0\,,\ \ \varepsilon_j\in \{-1,0,1\} \Longrightarrow \forall j\ \varepsilon_j=0\,.
\end{equation}
On se borne d'ordinaire aux produits (1) pour lesquels
\begin{equation}
|a_j|=a<1\,.
\end{equation}
C'est un outil essentiel pour l'\'etude des ensembles de Sidon $\Lambda$ $(\Lambda \subset \Z)$ d\'efinis par la condition que, pour les polyn\^omes trigonom\'etriques \`a spectre dans $\Lambda$, les normes dans $C(\T)$ et dans $A(\T)=\Ff \ell^1(\Z)$ sont \'equivalentes :
\begin{equation}
\sum_{\lambda\in \Lambda} |c_\lambda| \le S\ \sup_{t} \Big| \sum_{\lambda\in \Lambda}c_\lambda e^{i\lambda t}\Big|\,.
\end{equation}
La borne inf\'erieure des $S$ est la constante de Sidon $S(\Lambda)$. Le livre de Li et Queffelec \cite{liqu} contient une \'etude tr\`es compl\`ete des relations entre ensembles quasi--ind\'ependants et ensembles de Sidon, \`a partir des travaux de Pisier \cite{pis} et de Bourgain \cite{bou} ; j'en donne une version un peu all\'eg\'ee dans  arXiv : math. 0709.4386.

\vskip2mm

Voici les principaux r\'esultats :
\begin{eqnarray}
&&\textrm{Toute r\'eunion finie d'ensembles quasi--ind\'ependants}\\
\nonumber&&\textrm{est un ensemble de Sidon}\,.
\end{eqnarray}
\begin{eqnarray}
&&\textrm{(caract\'erisation de Bourgain) Une condition n\'ecessaire et suffisante}\\
\nonumber&&\textrm{pour qu'une partie } \Lambda  \textrm{ de } \Z \textrm{ soit un ensemble de Sidon est qu'il existe} \\
\nonumber&&\textrm{un } c=c(\Lambda)>0 \textrm{ tel que, pour toute mesure positive }\varpi \textrm{ sur } \Lambda, \textrm{ il existe}\\
\nonumber&& \textrm{une partie }\Lambda' \textrm{ de } \Lambda, \textrm{ quasi--ind\'ependante, telle que }\varpi(\Lambda')\ge c\, \varpi(\Lambda)\,.
\end{eqnarray}

La preuve de (46) utilise seulement les produits de Riesz, et elle donne une estimation de la constante de Sidon : si $\Lambda$ est la r\'eunion de $k$ ensembles quasi--ind\'ependants,
\begin{equation}
S(\Lambda) \le 3 \sqrt{3}\ k\  \sqrt{2k-1}\,.
\end{equation}
Pour $k=1$, cela donne $S(\Lambda) \le 5,20$, et on a l'estimation un peu meilleure $S(\Lambda)\le 4,27$ ; la meilleure borne n'est pas connue.

On ne sait si (46) est une caract\'erisation des ensembles de Sidon. La condition de Bourgain (47) est la meilleure approche actuelle de la r\'eciproque. Elle est bas\'ee sur la condition de Pisier, elle aussi n\'ecessaire et suffisante, dans laquelle $\varpi$ est la mesure de d\'ecompte sur une partie finie arbitraire de $\Lambda$. L\`a encore, les produits de Riesz sont constamment utilis\'es.

Quittons pour finir les produits de Riesz -- c'est la digression annonc\'ee -- pour examiner une relation entre ensembles quasi--ind\'ependants, ensembles de Sidon, et condition de maille. Une maille est un ensemble de la forme
\begin{equation}
M=M(k,\{\gamma_j\},E) = \Big\{\sum_{j=1}^k n_j \gamma_j\}\,,
\end{equation}
o\`u les $\gamma_j$ sont des entiers fix\'es $(\gamma_j\in \Z,\ j=1,2,\ldots k)$, et les $n_j$ des entiers tels que $(n_1,n_2,\ldots n_k)\in E$, ensemble fix\'e de $\Z^k$. On a su d\`es l'introduction des ensembles de Sidon qu'\`a chaque ensemble de Sidon $\Lambda$ est attach\'e un nombre positif $C$, ne d\'ependant que de $S(\Lambda)$, tel que pour toute maille $M$ de la forme (49) on ait
\begin{equation}
|\Lambda \cap M| \le C\, k\, \log (1+ \sup_{(n_1,n_2,\ldots n_k)\in E} (|n_1|+|n_2|+\cdots+|n_k|)\,.
\end{equation}
En particulier, si $E=\{-1,0,1\}^k$, on a
\begin{equation}
|\Lambda \cap M| \le C\, k\, \log(1+k)\,.
\end{equation}
Cela vaut en particulier si $\Lambda$ est quasi--ind\'ependant.

Pour montrer que (51) est inam\'eliorable, je me propose de construire un ensemble $\Lambda$ quasi--ind\'ependant tel que, pour tout entier $k$, il existe une maille $M$ de la forme (49) avec $E=\{-1,0,1\}^k$, telle~que
\begin{equation}
|\Lambda \cap M| \ge \frac{1}{4} \, k \, \log_2 k\,.
\end{equation}

La d\'efinition de la quasi--ind\'ependance s'\'etend imm\'ediatement  \`a des sous--ensembles de $\Z^n$. Quand $n$ est une puissance de 2, $n=2^\nu$, on va construire dans la maille $\{-1,0,1\}^n$ de $\Z^n$ un ensemble quasi--ind\'ependant assez riche. Commen\c cons par $n=2$. Les vecteurs colonnes de la matrice
\[
\left(
\begin{array}{ccc}
1  &1	   &1   \\
1  &-1   &0   \\
\end{array}
\right)
\]
sont q.i. (quasi--ind\'ependants) ; v\'erifions le en d\'etail. En effet, si
\begin{equation*}
\varepsilon_1 \begin{pmatrix}
1\\1
\end{pmatrix} +
\varepsilon_2 \begin{pmatrix}
1\\-1 
\end{pmatrix} +
\varepsilon_3 \begin{pmatrix}
1\\0
\end{pmatrix}  =
 \begin{pmatrix}
0\\0
\end{pmatrix} 
\end{equation*}
avec $\varepsilon_j\in \{-1,0,1\}$, on a d'abord $\varepsilon_1=\varepsilon_2$ (seconde ligne), puis $\varepsilon_3=0$ modulo 2 donc $\varepsilon_3=0$ (premi\`ere ligne), puis $\varepsilon_1=\varepsilon_2=0$ (ind\'ependance de $\begin{pmatrix}
1\\1
\end{pmatrix}$ et $\begin{pmatrix}
1\\-1
\end{pmatrix}$).

Lorsque $n=2^\nu$, on va construire par r\'ecurrence des matrices $A_\nu$ \`a $2^\nu$ lignes et $N_\nu$ colonnes, dont les colonnes sont dans $\{-1,0,1\}^n$ et sont q.i.\,. Pour $\nu=1$, c'est fait, avec $N_1=3$. On passe de $A_\nu$ \`a $A_{\nu+1}$ par le proc\'ed\'e figur\'e
\begin{equation*}
A_{\nu+1} = \begin{pmatrix}
A_\nu &A_\nu &I_\nu\\
A_\nu &-A_\nu &0
\end{pmatrix}
\end{equation*}
o\`u $I_\nu$ est la matrice unit\'e $2^\nu \times 2^\nu$. Les colonnes de $A_{\nu+1}$ sont dans $\{-1,0,1\}^{2n}$. Montrons qu'elles sont q.i. lorsque celles de $A_\nu$ le sont. Une relation lin\'eaire \`a coefficients $-1,0$ ou 1 entre les colonnes de $A_{\nu+1}$ s'\'ecrit, en posant $n=2^\nu$ et $N=N_\nu$,
\begin{eqnarray*}
&&A_\nu(\varepsilon_1^1, \varepsilon_1^2, \ldots \varepsilon_1^N)^t +A_\nu(\varepsilon_2^1,\varepsilon_2^2,\ldots \varepsilon_2^N)^t + I_\nu(\varepsilon_3^1,\varepsilon_3^2,\ldots \varepsilon_3^n)^t =0\,,\\
&&A_\nu(\varepsilon_1^1, \varepsilon_1^2, \ldots \varepsilon_1^N)^t -A_\nu(\varepsilon_2^1,\varepsilon_2^2,\ldots \varepsilon_2^N)^t  =0\,,
\end{eqnarray*}
les $\varepsilon$ valant $-1,0$ ou 1. En ajoutant, on voit que les lignes de $I_\nu(\varepsilon_3^1,\varepsilon_3^2,\ldots \varepsilon_3^n)^t$ sont nulles modulo 2, donc nulles, donc $\varepsilon_3^1 =\varepsilon_3^2=\cdots =\varepsilon_3^n =0$, donc
\begin{equation*}
A_\nu(\varepsilon_1^1,\varepsilon_1^2,\ldots \varepsilon_1^N)^t = A_\nu(\varepsilon_2^1,\varepsilon_2^2,\ldots \varepsilon_2^N)=0\,,
\end{equation*}
et la quasi--ind\'ependance des colonnes de $A_\nu$ entra\^{i}ne que tous les $\varepsilon$ sont nuls. Les colonnes de $A_{\nu+1}$ sont donc bien q.i.\,.

Calculons $N_\nu$. Partons de $N_0=1$. La construction donne
\begin{equation*}
N_\nu = 2N_{\nu-1} + 2^{\nu-1}\,,
\end{equation*}
soit
\begin{equation*}
2^{-\nu} N_\nu = 2^{-(\nu-1)} N_{\nu-1} + \frac{1}{2} =\cdots = N_0 + \frac{\nu}{2}
\end{equation*}
donc
\begin{equation*}
N_\nu = 2^{\nu-1}(2+\nu)
\end{equation*}

Pour d\'efinir $\Lambda$, on choisit dans $\Z$ une suite $(\beta_j)_{j\ge 1}$ tr\`es dissoci\'ee dans le sens suivant : il n'y a pas de relation lin\'eaire non triviale du type $\sum n_j \beta_j=0$ (somme finie) avec $n_j\in \Z$ et $|n_j|\le N_\nu$ quand $2^\nu \le j < 2^{\nu+1}$ ($\nu\ge 1$). On peut construire une telle suite par r\'ecurrence sur $\nu$. Pour
\begin{equation*}
\sum_{i=1}^{\nu-1} N_i < \ell \le \sum_{i=1}^\nu N_i
\end{equation*}
on d\'efinit le vecteur ligne $(\gamma_\ell)$ comme
\begin{equation*}
(\gamma_\ell)= (\beta_{2^\nu}, \beta_{2^\nu+1}, \ldots \beta_{2^{\nu+1}-1})A_\nu
\end{equation*}
(c'est $A_\nu$ \'etal\'e sur $\Z$ au moyen des $\beta_j$). La suite cherch\'ee $\Lambda$ est la suite $\gamma_\ell$ $(\ell\ge 1)$. Elle est q.i. parce que toute expression de la forme $\sum \varepsilon_\ell \gamma_\ell$ $(\varepsilon_\ell \in \{-1,0,1\})$ s'\'ecrit $\sum n_j\beta_j$ avec $|n_j|\le N_\nu$ quand $2^\nu \le j <2^{\nu+1}$. Elle a $N_\nu$ termes dans la maille
\begin{equation*}
M=\{\Sigma \varepsilon_j \beta_j\,;\ \varepsilon_j \in \{-1,0,1\}\,;\ 2^\nu \le j<2^{\nu+1}\}\,,
\end{equation*}
qui est bien de la forme (49) avec $k=2^\nu$, $\gamma_i = \beta_{i+2^\nu+1}$, et $E=\{-1,0,1\}^k$. Pour une telle maille $M$, on a
\begin{equation*}
|\Lambda\cap M| = N_\nu = \frac{1}{2} \, 2^\nu (2+\nu) > \frac{1}{2} \ k\ \log_2 k\,.
\end{equation*}
Mais $M$ est aussi une maille de la forme (49) avec $k$ quelconque entre $2^\nu$ et $2^{\nu+1}$. Ainsi, pour tout entier $k$, on peut trouver $M$ de la forme (49) telle que
\begin{equation*}
|\Lambda \cap M| > \frac{1}{4}\, k\ \log_2 k\,,
\end{equation*}
ce qu'on voulait montrer.

\vskip2mm

Le cadre naturel pour les ensembles q.i. et pour les ensembles de Sidon est celui des groupes ab\'eliens discrets infinis. Si $\Gamma$ est un tel groupe, on transcrit imm\'ediatement pour $\Gamma$ les d\'efinitions donn\'ees pour $\Z$. Les propositions (46) et (47) restent valables. La derni\`ere construction reste valable si $\Gamma$ contient des \'el\'ements d'ordre arbitrairement grand, et seulement dans ce cas.

Lorsque $\Gamma$ est le groupe dual de $\T^\N$, Pisier a obtenu le m\^eme r\'esultat que notre construction par une m\'ethode non--constructive (\cite{pis}, proposition~7.3). Le corollaire~3.3 de \cite{pis} montre que le r\'esultat est en d\'efaut quand $\Gamma$ est le groupe dual de $\prod\limits_{j=1}^\infty(\Z/p_j\Z)$, o\`u $p_j$ est une suite born\'ee d'entiers ; dans ce cas, tout ensemble de Sidon rencontre une maille de la forme (49) en au plus $Ck$ points, $C$ ne d\'ependant que de l'ensemble de Sidon.

La construction et les remarques sont consultables sur arXiv : math. 0709.4383.

\vskip3mm

Il me reste \`a faire un commentaire sur les commentaires. Tous essayent de mettre en valeur le r\^ole de pionnier de Peyri\`ere, soit dans les \'enonc\'es, soit dans les m\'ethodes. Le commentaire sur $C$ est accompagn\'e d'une tentative, peu concluante, pour rattacher les analyses dimensionnelle et multifractale d'un produit de Riesz \`a sa s\'erie de Fourier. Une digression s'ajoute au commentaire sur $D$, \`a savoir la construction d'un ensemble quasi--ind\'ependant aux propri\'et\'es tr\`es diff\'erentes de celles de $\{2^j\}$; elle se rattache naturellement \`a $D$ parce que les r\'esultats de Peyri\`ere sur $\{2^j\}$ sont l'avant--garde de la th\'eorie des ensembles quasi--ind\'ependants.

Une derni\`ere remarque s'impose : les travaux de Peyri\`ere sur les produits de Riesz ne se limitent pas \`a ses articles de 1975 \`a 1990 cit\'es en r\'ef\'erences. Les id\'ees et les r\'esultats les plus saillants se trouvent d\'ej\`a dans la note aux Comptes Rendus \og Sur les produits de Riesz\fg\  de 1973 (volume 276, pp.~1453--1455), et dans l'expos\'e 22 du S\'eminaire de Th\'eorie des Nombres de l'Universit\'e de Bordeaux 1 de 1973-1974, intitul\'e \og Utilisation des produits de Riesz pour minorer la dimension de Hausdorff de certains ensembles\fg.

\eject

\end{document}